\documentclass[11pt]{amsart}
\usepackage{amsmath,amssymb,amsthm,mathtools}
\usepackage{aliascnt}
\usepackage[margin=1in]{geometry}
\usepackage{microtype}
\usepackage{float}
\usepackage{graphicx}
\usepackage{hyperref}
\usepackage[noabbrev,capitalize]{cleveref}

\newtheorem{theorem}{Theorem}

\newaliascnt{lemma}{theorem}
\newtheorem{lemma}[lemma]{Lemma}
\aliascntresetthe{lemma}
\crefname{lemma}{Lemma}{Lemmas}
\Crefname{lemma}{Lemma}{Lemmas}

\newaliascnt{proposition}{theorem}
\newtheorem{proposition}[proposition]{Proposition}
\aliascntresetthe{proposition}
\crefname{proposition}{Proposition}{Propositions}
\Crefname{proposition}{Proposition}{Propositions}

\newaliascnt{corollary}{theorem}
\newtheorem{corollary}[corollary]{Corollary}
\aliascntresetthe{corollary}
\crefname{corollary}{Corollary}{Corollaries}
\Crefname{corollary}{Corollary}{Corollaries}
\theoremstyle{remark}
\newtheorem*{remark}{Remark}

\numberwithin{equation}{section}

\newcommand{\C}{\mathbb C}
\newcommand{\F}{\mathbb F}
\newcommand{\Q}{\mathbb Q}
\newcommand{\Z}{\mathbb Z}
\newcommand{\cD}{\mathcal D}
\newcommand{\cO}{\mathcal O}
\newcommand{\cQ}{\mathcal Q}
\newcommand{\cX}{\mathcal X}
\newcommand{\Cl}{\operatorname{Cl}}
\newcommand{\End}{\operatorname{End}}
\newcommand{\Res}{\operatorname{Res}}
\newcommand{\leg}[2]{\left(\frac{#1}{#2}\right)}

\newcommand{\appendixcredit}{\normalfont\normalsize (Appendix by Qi-Yang Zheng)}

\title[Parity in quadratic progressions]
{Parity of the partition function in quadratic progressions\\
\texorpdfstring{\appendixcredit}{(Appendix by Qi-Yang Zheng)}}
\dedicatory{In memory of Dick Gross}


\author{Ken Ono}
\address{Axiom Math, 124 University Avenue, Palo Alto, CA 94301, USA}
\email{ken@axiommath.ai}
\address{Department of Mathematics, University of Virginia, Charlottesville, VA 22904, USA}
\email{ko5wk@virginia.edu}

\author{Ashvin Swaminathan}
\address{Axiom Math, 124 University Avenue, Palo Alto, CA 94301, USA}
\email{ashvin@axiommath.ai}

\begin{document}
\begin{abstract}
The parity of the partition function $p(n)$ is one of the most
stubborn problems in number theory.  In 2010, the first author conjectured, for square-free $1<D\equiv 23\pmod{24}$, that the values
$p\!\left(\frac{Dm^2+1}{24}\right)$,
as $m$ ranges over positive integers with $(m,6)=1$, include infinitely
many even and infinitely many odd terms.  We prove
this conjecture.
The key new idea is geometric.  Logarithmic derivatives of
\emph{twisted Borcherds products} built from Ramanujan's third-order
mock theta functions recast the problem in terms of CM points of
discriminant $-D$ on $X_0(6)$.  The uniqueness of canonical lifts from
characteristic $2$ to characteristic zero shows that the CM points
supporting the poles remain distinct after reduction modulo $2$.
This fact, combined with an Eisenstein series comparison and a Galois
representation argument, rules out constant parity and gives infinitely many
values of each parity.
The method applies to the coefficients of analogous generalized
twisted Borcherds products. 
These results imply that
\[\#\{0\leq n\leq X: p(n)\text{ is odd}\}\gg\sqrt X,\]
which is now the best known lower bound for odd values of $p(n)$. The algebraic identities at the heart of this paper, as well as the improved lower bound for odd partition numbers, have been formalized in Lean by AxiomProver.
\end{abstract}
\maketitle

\section{Introduction and statement of results}

A \emph{partition} of a positive integer $n$ is a way to write $n$ as a sum
of positive integers, disregarding order.  The number of such ways is the
partition function $p(n)$.  For example, $p(4)=5$.  In 1919, Ramanujan proved
\cite{Ramanujan} the first two of the three striking congruences
\[
 p(5n+4)\equiv0\pmod5,\qquad
 p(7n+5)\equiv0\pmod7,\qquad
 p(11n+6)\equiv0\pmod{11},
\]
and he conjectured the third.  They are glimpses of families of congruences
modulo powers of $5$, $7$, and $11$, revealing unexpected regularities in the coefficients of
\begin{equation}\label{eq:Pgen}
 P(q):=\sum_{n\geq0}p(n)q^n
 =\prod_{n=1}^{\infty}\frac{1}{1-q^n}.
\end{equation}
These congruences launched a far-reaching theory linking partitions and
modular forms 
\cite{AhlgrenOno2001,AndrewsTP,Atkin,Atkin1968,Ono2000,Watson1938}.

Despite a century of progress, our understanding of the parity of $p(n)$
remains stubbornly limited.  The Parkin--Shanks conjecture predicts that
$p(n)$ is even and odd with natural density $1/2$ each
\cite{ParkinShanks1967}.  Prior to this paper, the best unconditional bounds,
building on work of Nicolas and Nicolas--Serre, were due to
Bella\"iche--Green--Soundararajan and Bella\"iche--Nicolas.   As $X\rightarrow +\infty$, their results (see \cite{BGS2018,BN,Nicolas,NicolasIJNT,NicolasSerre}) give
\begin{displaymath}
\begin{split}
 \#\{0\leq n\leq X:p(n)\text{ is odd}\}
 & \gg \frac{\sqrt X}{\log\log X},\\
 \#\{0\leq n\leq X:p(n)\text{ is even}\}
 &\gg \sqrt X\cdot \log\log X.
\end{split}
\end{displaymath}

A natural way to probe this apparent randomness is to restrict the argument
of $p(n)$ to the values of a polynomial.  The first author has conjectured more
broadly that if
$g(x)$ is a nonconstant integer-valued polynomial that is eventually positive,
then $p(g(n))$ is even for infinitely many $n$ and odd for infinitely many
$n$ (see \cite{Ono2010}).
For polynomials of degree $1$, this is precisely Subbarao's conjecture: every
arithmetic progression contains infinitely many even and infinitely many odd
values of $p(n)$.  In \cite{Ono1996}, the first author proved the even case in full
and showed that one odd value in a fixed progression forces infinitely many
odd values there.  Radu \cite{Radu2012} later proved that every arithmetic
progression contains such an odd value, thereby completing the conjecture.

For polynomials of degree at least $2$, essentially nothing is known.  In
\cite{Ono2010}, the first author singled out the special quadratic progressions
\[
 n=\frac{Dm^2+1}{24},\qquad (m,6)=1,
\]
where $D\equiv23\pmod{24}$ is square-free, and conjectured that both
parities occur infinitely often on each of these progressions.  Partial
progress toward this conjecture was made in earlier work
\cite{Ono2010,Ono2024}.
Note that $(m,6)=1$ forces
$m^2\equiv1\pmod{24}$, so that $Dm^2+1\equiv0\pmod{24}$ and these progressions
consist of integers.  We call a positive integer $m$ \emph{admissible} when
$(m,6)=1$.  

We solve the conjecture outright, using new ideas from geometry.  The first
theorem states the result together with effective bounds for the first
occurrence of each parity.

\begin{theorem}\label{thm:main}
 If $1<D\equiv23\pmod{24}$ is square-free, then both parities occur
 infinitely often among
 \[
  p\!\left(\frac{Dm^2+1}{24}\right),
  \]
  where $(m,6)=1$.
 Furthermore, if $h(-D)$ is the class number of
 $\Q(\sqrt{-D})$, and if $m_{\mathrm{odd}}$ and $m_{\mathrm{even}}$ denote
 the least positive admissible integers producing an odd and an even value,
 respectively, then
 \begin{equation}\label{eq:first-bounds}
   m_{\mathrm{odd}}\leq 12h(-D)+2 \qquad\text{and}\qquad
   m_{\mathrm{even}} \leq
   \bigl(12h(-D)+2\bigr)\prod_{\ell\mid D}(\ell+1).
 \end{equation}
\end{theorem}

Averaging the bound for $m_{\mathrm{odd}}$ in
\eqref{eq:first-bounds} over $D$ produces many odd values of
$p(n)$, enough to give the first lower bound of order $\sqrt X$ for their count.

\begin{corollary}\label{cor:sqrt}
 If we let $N_{\mathrm{odd}}(X):=\#\{0\leq n\leq X: p(n)\text{ is odd}\}$, then we have
 \[
  \liminf_{X\to\infty}\frac{N_{\mathrm{odd}}(X)}{\sqrt X}
  \ \geq\ \frac{243}{64\sqrt6\,\pi^5}.
 \]
 In particular, we have
 \[ \# \left \{ 0\leq n\leq X \ : \ p(n) \ {\text {is odd}}\right \}\gg\sqrt X.
 \]
\end{corollary}

\noindent This corollary is due to Qi-Yang Zheng, and his proof is given in  Appendix \ref{sec:appendix}.  It
improves the bound of Bella\"iche--Green--Soundararajan quoted above by a
factor of $\log\log X$.

At first sight, the partition values in \cref{thm:main} seem far removed
from geometry.  The bridge begins with an elementary identity: modulo $4$,
the generating function for $p(n)$ agrees with Ramanujan's third-order mock
theta function $f(q)$, so the parities in \cref{thm:main} are encoded by
selected coefficients of $f$.  This mock theta function is the holomorphic
part of a weight $1/2$ harmonic Maass form.  For each $D$, the twisted Borcherds lift of
Bruinier and the first author \cite{BruinierOnoAnnals} gives a
meromorphic product $\Psi_D$ of weight $0$ for $\Gamma_0(6)$.  The infinite
product remembers the relevant partition values, while its divisor consists of two packets
supported at CM points of discriminant $-D$.  Its normalized logarithmic
differential therefore has two complementary descriptions: its Fourier
expansion is a Lambert series encoding the desired parities, and its poles
are governed by CM geometry.

The new geometric ingredient is the behavior of these poles in characteristic $2$.
The two CM packets in the divisor have a common $2$-orientation, and canonical-lift
rigidity shows that distinct CM points remain distinct after reduction at
the compatible prime above $2$.  The poles at these points are simple, and
their residues are $2$-adic units.  Therefore, the poles survive modulo
$2$.  This rules out the
possibility that all the partition values are even.  If they were all odd,
the same Lambert series would be congruent to a holomorphic Eisenstein
series, which is impossible because the former differential retains a CM
pole.  Section~\ref{sec:geometry} develops the geometry behind this argument.

The same reasoning gives a nonconstancy criterion for sequences arising
from generalized twisted Borcherds products.  We state it in geometric form
because integrality of a holomorphic part and the existence of a Borcherds
product alone do not guarantee the integral and characteristic $2$ properties
needed in the proof.

Let $N=2N_0$ with $N_0$ odd, and let $C_\infty$ be the ordinary component of
the Deligne--Rapoport special fiber of $X_0(N)$ containing $\infty$.  For an
odd square-free $D$ coprime to $N$, let $B_D(m)\in\Z$ be defined for
$(m,N)=1$, and put
\[
 \mathcal L_D(B;q):=
 \sum_{\substack{m\geq1\\(m,N)=1}}B_D(m)
 \sum_{\substack{n\geq1\\(n,D)=1}}q^{mn}.
\]

The following theorem isolates the precise inputs from geometry and modular
forms used in the proof.  Its conclusion is deliberately one of
nonconstancy; an additional arithmetic propagation result is needed to turn
one occurrence of each parity into infinitely many.

\begin{theorem}\label{thm:general}
 Assume that \(-D\) is a fundamental discriminant and that \(2\) splits
in \(K=\Q(\sqrt{-D})\).
 \begin{enumerate}
 \item[(G1)] There is a meromorphic differential $\omega_D$ on
 $X_0(N)_{\Q_2}$ whose expansion at $\infty$ is integral with vanishing
 constant term (so that $\omega_D$ is regular at the cusp $\infty$),
 such that
 \[
   \omega_D\equiv \mathcal L_D(B;q)\cdot \frac{dq}{q}\pmod2.
 \]
 \item[(G2)] The noncuspidal polar divisor of $\omega_D$ is nonempty and
 consists of simple maximal-order CM points of discriminant $-D$.  Their
 residues are $2$-adic units, their prime-to-$2$ level structures are
 CM-stable, and their order $2$ subgroups are annihilated by one fixed prime
 $\mathfrak p_2\mid2$ of $K$.
 \item[(G3)] There is a holomorphic weight $2$ modular form
 $\mathcal E_D$ on $\Gamma_0(ND)$, defined over $\Q_2$ and integral at
 $\infty$, such that
 \[
  \mathcal E_D(q)\equiv
  \sum_{\substack{m\geq1\\(m,N)=1}}
  \sum_{\substack{n\geq1\\(n,D)=1}}q^{mn}\pmod2.
 \]
 \end{enumerate}
 Choose the $2$-adic CM place for which the subgroup in \emph{(G2)} is
 connected, and suppose all the CM poles, their CM actions, and their level
 structures in \emph{(G2)} are defined over a common finite unramified
 extension of $\Q_2$.  Then the sequence $B_D(m)$
 is not constant modulo $2$.
 If, in addition, the occurrence of either parity at one index $m$ with
 $(m,N)=1$ is known to imply its occurrence at infinitely many such indices,
 then both parities occur infinitely often.
\end{theorem}

\begin{remark}
The word \emph{nonconstancy} in \cref{thm:general} is essential.
The geometric argument rules out the possibility that all the $B_D(m)$ have
the same parity:
\[
 B_D(m)\equiv 0\pmod 2
 \qquad\text{or}\qquad
 B_D(m)\equiv 1\pmod 2,
\]
but it gives no relation between the values of $B_D(m)$ at distinct
indices.  Thus, in principle, hypotheses \emph{(G1)--(G3)} are
compatible with a sequence having only finitely many terms of one
parity.
For $p(n)$, the passage from one occurrence to
infinitely many is supplied by a separate arithmetic argument
\cite[Theorem~1.2 and Section~3.4]{Ono2010}.  After clearing the CM poles
to obtain holomorphic modular forms modulo $2$, that argument uses the
local nilpotence of certain Hecke algebras, together with distribution
results derived from Galois representations and the Chebotarev density
theorem.  The required nilpotence is a special small-level phenomenon,
closely related to the absence of certain nontrivial mod $2$ Galois
representations, and need not hold at general levels.  Consequently,
an infinitude theorem for general Borcherds data requires an additional
propagation theorem, typically involving Hecke operators or Galois
representations, and is not a formal consequence of \emph{(G1)--(G3)}.
\end{remark}

The paper is organized as follows.  Section~\ref{sec:mock} recalls the
relevant theory of generalized twisted Borcherds products, gives the exact
two-packet divisor, and proves the partition Lambert series identity.
Section~\ref{sec:geometry} develops the ordinary characteristic $2$ geometry,
including CM orientation, separation, and survival of poles.
Section~\ref{sec:main-proof} proves
\cref{thm:main}.  Section~\ref{sec:general-proof} sketches the proof of
\cref{thm:general}.  \Cref{sec:appendix} contains the proof of \cref{cor:sqrt}
due to Zheng, and \cref{sec:lean} describes the Lean formalizations of the
algebraic identities required for \cref{thm:main}, as well as of the proof of
\cref{cor:sqrt}.

\section*{Acknowledgements}
\noindent The authors thank George Andrews and Jan Hendrik Bruinier for their comments on earlier versions of this paper. They also thank Kenny Lau for his assistance with the AxiomProver. The first author thanks the Thomas Jefferson Fund,
the NSF (DMS-2002265 and DMS-2055118), and the Simons Foundation
(SFI-MPS-TSM-00013279) for their generous support.

\section{Mock theta functions, Borcherds products, and
\texorpdfstring{$p(n)$}{p(n)}}
\label{sec:mock}

This section assembles the ingredients from \cite{BruinierOnoAnnals,Ono2010}
that connect the partition function to twisted Heegner divisors.  The
combinatorial identity and the components of the Maass form are recalled from
those sources.  The divisor computation in \S\ref{subsec:divisor} is carried
out in full, because the two-packet structure and the unit multiplicities in
\eqref{eq:exact-divisor} are used in an essential way in
Section~\ref{sec:geometry}.  The central
combinatorial observation is that Ramanujan's third-order mock theta function
encodes $p(n)$ modulo $2$.  The products used below belong to the generalized
Borcherds theory originating in \cite{Borcherds1,Borcherds2}.
Throughout, $\tau\in\mathbb H$, $q=e^{2\pi i\tau}$, and
$(a;q)_n:=\prod_{j=0}^{n-1}(1-aq^j)$.  We will not need the general
construction of the lift.  The explicit components of the Maass form, the
product, and its divisor stated below contain all of the input used later.

\subsection{Combinatorial considerations}

A partition can be represented by its Ferrers diagram.  The largest square in
the upper-left corner is its Durfee square; its boundary separates the
partition into a square and two partitions whose parts do not exceed the side
length.  For example, the partition $5+4+2+1$ has Ferrers diagram

\begin{figure}[H]
  \centering
    \includegraphics[height=35mm]{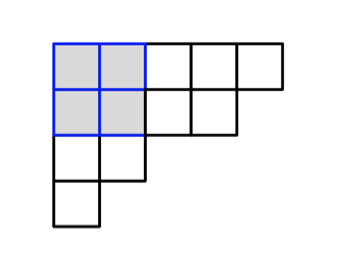}
  \caption{A partition with its Durfee square.}
\end{figure}
This partition decomposes into a $2\times2$ Durfee square and the two
partitions $2+2+1$ and $2+1$, where the first is obtained by transposing the
portion to the right of the Durfee square.

Ramanujan's third-order mock theta function is
\begin{equation}\label{eq:fq}
 f(q):=1+\sum_{n=1}^{\infty}
 \frac{q^{n^2}}{(1+q)^2(1+q^2)^2\cdots(1+q^n)^2}.
\end{equation}

The Durfee-square decomposition gives the elementary congruence that begins
the passage from partitions to mock modular forms.

\begin{lemma}\label{lem:durfee}
 As formal power series, we have
 \[
  P(q)=1+\sum_{m=1}^{\infty}
  \frac{q^{m^2}}{(1-q)^2(1-q^2)^2\cdots(1-q^m)^2}.
 \]
 In particular, $P(q)\equiv f(q)\pmod4$.
\end{lemma}

\begin{proof}
 The series
 \[
  \frac{1}{(1-q)(1-q^2)\cdots(1-q^m)}
 \]
 generates partitions whose parts do not exceed $m$.  Squaring it and
 multiplying by $q^{m^2}$ records the two partitions adjacent to a Durfee
 square of side $m$.  Summing over $m$ proves the identity.  Finally,
 $(1-q^j)^2\equiv(1+q^j)^2\pmod4$, and power series with constant term $1$
 that agree modulo $4$ have reciprocals that agree modulo $4$; comparing with
 \eqref{eq:fq} gives the congruence.
\end{proof}

\subsection{The mock-theta harmonic Maass form}

A harmonic weak Maass form is a real-analytic modular form with controlled
growth at the cusps that is annihilated by the appropriate
weight Laplacian.  It has a canonical decomposition into a holomorphic part
$H^+$, given by a Fourier series, and a nonholomorphic part.  In the
vector-valued setting its components transform together through a Weil
representation.  This framework is useful here because the coefficients of
$H^+$ become exponents in a Borcherds product, while the negative-index
terms determine the divisor of that product.  See
\cite[Section~3]{BruinierFunke} for harmonic weak Maass forms and
\cite[Sections~6 and~8.2]{BruinierOnoAnnals} for the lift used below.

Along with $f(q)$, consider
\begin{equation}\label{eq:omega-mock}
 \omega(q):=\sum_{n=0}^{\infty}
 \frac{q^{2n^2+2n}}{(q;q^2)_{n+1}^2}.
\end{equation}
Both $f$ and $\omega$ have integral coefficients.  For $0\leq j\leq11$,
let
\[
 H_j^+(\tau)=\sum_{n\gg-\infty}C(j;n)q^n
\]
be given by
\begin{equation}\label{eq:weil}
H_j^+(\tau)=
\begin{cases}
0,&j=0,3,6,9,\\
q^{-1}f(q^{24}),&j=1,7,\\
-q^{-1}f(q^{24}),&j=5,11,\\
2q^8\bigl(-\omega(q^{12})+\omega(-q^{12})\bigr),&j=2,\\
-2q^8\bigl(\omega(q^{12})+\omega(-q^{12})\bigr),&j=4,\\
2q^8\bigl(\omega(q^{12})+\omega(-q^{12})\bigr),&j=8,\\
2q^8\bigl(\omega(q^{12})-\omega(-q^{12})\bigr),&j=10.
\end{cases}
\end{equation}
After adjoining the nonholomorphic period integrals described in
\cite[Section~8.2]{BruinierOnoAnnals}, these components form a vector-valued
weight $1/2$ harmonic weak Maass form for the Weil representation attached to
$X_0(6)$.

\subsection{The exact Borcherds divisor}\label{subsec:divisor}

We now pass from the mock-theta coefficients to geometry.  Throughout this
subsection, assume that $D>1$ is square-free with
$D\equiv23\pmod{24}$, put $K:=\Q(\sqrt{-D})$, and use the
parameter $r=1$, which is permitted because $-D\equiv1\pmod{24}$.  Define
\begin{equation}\label{eq:PD}
 P_D(X):=\prod_{b\bmod D}
 \bigl(1-e(-b/D)X\bigr)^{\leg{-D}{b}},
 \qquad e(x):=e^{2\pi i x},
\end{equation}
and
\begin{equation}\label{eq:PsiD}
 \Psi_D(z):=\prod_{m=1}^{\infty}
 P_D(q^m)^{C(\overline m;Dm^2)},
\end{equation}
where $\overline m$ is the residue class of $m$ modulo $12$.

We now describe its divisor in elementary terms.  For $s\in\{1,5\}$, let
$\cQ_{D,s}$ be the set of $\Gamma_0(6)$-classes of primitive positive-definite
forms
\[
 Q=[a,b,c],\qquad b^2-4ac=-D,\qquad 6\mid a,\qquad b\equiv s\pmod{12},
\]
and put
$P_Q:=\frac{-b+\sqrt{-D}}{2a}\in\mathbb H.$
In the moduli interpretation, $P_Q$ is the point of $X_0(6)$ corresponding to
the enhanced elliptic curve
\[
 \bigl(E_Q,\ \langle\tfrac16\rangle\bigr),
 \qquad E_Q:=\C/(\Z\tau_Q+\Z),\quad \tau_Q=P_Q,
\]
and $E_Q$ has CM by the maximal order $\cO_K$ because $Q$ is primitive and
$-D$ is fundamental.  Since $6\mid a$, the residue $b\pmod{12}$ is an
invariant of the $\Gamma_0(6)$-class of $Q$ (indeed $b\bmod 2N$ is a
$\Gamma_0(N)$-class invariant of forms with $N\mid a$), so the two packets
$\cQ_{D,1}$ and $\cQ_{D,5}$ are disjoint sets of points on $X_0(6)$.
Each packet is a torsor for $\Cl(\cO_K)$ and therefore contains $h(-D)$
points (for example, see \cite{GrossX0N}).  The \emph{generalized genus character}
in the twisted Heegner divisor
supplies a sign $\epsilon_D(Q):=\chi_{-D}(Q)\in\{\pm1\}$ for every such
primitive form.  In fact $\epsilon_D\equiv+1$ here.  The twisting discriminant is $-D$ and
the forms occurring have discriminant $-D$, so the relevant factorization
in the sense of \cite[(4.1)]{BruinierOnoAnnals} is the trivial one
$-D=(-D)\cdot1$.  Equivalently, $\chi_{-D}(Q)=\left(\frac{-D}{n}\right)$
for any $n$ prime to $D$ represented by $Q$; since $Q$ has discriminant
$-D$, any such $n$ satisfies $\left(\frac{-D}{n}\right)=+1$.  We retain
the notation $\epsilon_D$ below, but the signs in \eqref{eq:exact-divisor}
are carried entirely by the packet labels.

The Borcherds product theorem identifies these two explicitly described
packets as the complete divisor of $\Psi_D$.  The assertion about unit
multiplicities will later become the assertion that the logarithmic
differential has unit residues.

\begin{theorem}[Bruinier--Ono]\label{thm:borcherds}
 The product $\Psi_D$ is a weight $0$ meromorphic modular form for
 $\Gamma_0(6)$ with unitary character.
 Its divisor has no cuspidal part and is
 \begin{equation}\label{eq:exact-divisor}
  \operatorname{div}(\Psi_D)
  =\sum_{Q\in\cQ_{D,1}}\epsilon_D(Q)P_Q
  -\sum_{Q\in\cQ_{D,5}}\epsilon_D(Q)P_Q.
 \end{equation}
 Every geometric multiplicity in \eqref{eq:exact-divisor} is $+1$ or $-1$.
Moreover, put $\Theta:=q\frac{d}{dq}$.  Then
 \begin{equation}\label{eq:FD}
  F_D(z):=\frac{1}{\sqrt{-D}}
  \frac{\Theta\Psi_D(z)}{\Psi_D(z)}
 \end{equation}
 is a weight $2$ meromorphic modular form on $\Gamma_0(6)$ with integral
 Fourier coefficients and simple poles at the points occurring in
 \eqref{eq:exact-divisor}.
\end{theorem}

\begin{proof}
 In the notation of \cite[Section~8.2]{BruinierOnoAnnals}, the product is the
 lift with twisting discriminant $-D$ and parameter $r=1$, and its divisor is
 \[
  2Z_{-D,1}\!\left(-\frac1{24},\frac1{12}\right)
  -2Z_{-D,1}\!\left(-\frac1{24},\frac5{12}\right).
 \]
 Under the standard identification of the lattice vectors in these two
 Heegner divisors with positive-definite forms $[a,b,c]$ satisfying $6\mid a$,
 the two summands are precisely the packets $b\equiv1\pmod {12}$ and
 $b\equiv5\pmod {12}$ used above.  This also gives the class-group torsor
 assertion and the cardinality of each packet.

 We next account for the multiplicities.  By
 \cite[(5.3)]{BruinierOnoAnnals}, the divisor is the sum of
 $c^+(m,h)Z_{-D,1}(m,h)$ over all $h$ in the discriminant group $L'/L$ of the
 lattice $L$ of \cite[Section~5]{BruinierOnoAnnals}, where the $c^+(m,h)$ are
 the coefficients of the holomorphic part of the input Maass form; for twisting
 discriminant $-D<0$ one has $c^+(m,h)=-c^+(m,-h)$, and by
 \cite[Lemma~5.1(iv)]{BruinierOnoAnnals} one has
 $Z_{-D,1}(m,-h)=\operatorname{sgn}(-D)Z_{-D,1}(m,h)=-Z_{-D,1}(m,h)$.
 By \eqref{eq:weil}, the four components with $j\in\{1,5,7,11\}$ pair as
 $\{1,11\}$ and $\{5,7\}$ under $h\mapsto-h$, and the two members of each
 pair therefore contribute equally.  This is the source of the factor $2$
 in the divisor displayed above.  
 
 In \cite[(5.1)]{BruinierOnoAnnals} each
 Heegner point is weighted by $\chi_{-D}(\lambda)/w(\lambda)$, where $\lambda$
 runs over the Heegner lattice vectors representing the points of the divisor
 and $w(\lambda)$ is the order of the stabilizer of $\lambda$ in $\Gamma_0(6)$.
 Here, in the notation of \cite[Section~5]{BruinierOnoAnnals}, we have
 $m=-\tfrac1{24}$ and $N=6$, so $d:=4Nm\operatorname{sgn}(-D)=1$; hence
 $d\cdot(-D)=-D<-4$ and $w(\lambda)=2$.
 The doubling cancels this weight exactly, so every $\lambda$ contributes
 $\pm1$.  Since the two packets are disjoint, as noted above, every point
 of $\cQ_{D,1}\cup\cQ_{D,5}$ therefore occurs with multiplicity $\pm1$
 (for example, see
 \cite[Section~8.2]{BruinierOnoAnnals}).  Since the forms are primitive and $-D$ is fundamental, the generalized genus character $\chi_{-D}$ is nonzero, and
 by the discussion preceding \cref{thm:borcherds} it equals $+1$ on every
 form occurring.  This gives
 \eqref{eq:exact-divisor} with unit geometric multiplicities. In
 particular, each packet contains $h(-D)$ points, so
 $\operatorname{div}(\Psi_D)$ has degree $h(-D)-h(-D)=0$, as it must for a
 weight $0$ meromorphic modular form.  Because the twisting discriminant is not $1$, all Weyl
 vectors vanish, so there is no cuspidal contribution to the divisor.
 
 The modularity and product expansion are Theorems~6.1--6.2 and
 Section~8.2 of \cite{BruinierOnoAnnals}.  Logarithmic differentiation removes
 the character and gives a weight $2$ meromorphic modular form.  Integrality
 of the Fourier coefficients of $F_D$ and the assertion about its simple
 poles are also recorded in \cite[Theorem~2.2 and equation~(2.6)]{Ono2010}.
 Integrality also follows directly from the exact expansion
 \eqref{eq:exact-FD-expansion} below, and the simplicity of the poles follows from the unit multiplicities in \eqref{eq:exact-divisor}: at a zero or pole of order $\pm1$, the logarithmic derivative has a simple pole with residue $\pm1$.  In particular, the
 $q$-expansion principle shows that $F_D$ is defined over $\Q$.
\end{proof}

The use of $F_D$, rather than a global reduction of $\Psi_D$, is important.
The logarithmic differential is insensitive to a multiplier and is the
object for which both the integral Fourier expansion and the local residues
are directly available.  The following corollary records the residues of
the logarithmic differential.

\begin{corollary}\label{cor:unit-residues}
 The differential
 \begin{equation}\label{eq:omegaD}
  \omega_D:=F_D(q)\frac{dq}{q}
  =\frac{1}{\sqrt{-D}}\,d\log\Psi_D
 \end{equation}
 has residue
 $\Res_{P_Q}(\omega_D)=\pm1/\sqrt{-D}$
 at every point in \eqref{eq:exact-divisor}.  In particular, every residue
 is a $2$-adic unit at every place of $\Q(\sqrt{-D})$ above $2$.
\end{corollary}

\begin{proof}
 In any local coordinate $t$ at $P_Q$ (the points $P_Q$ have trivial
 stabilizer in $\Gamma_0(6)/\{\pm I\}$ because $-D<-4$, so $t$ is an honest
 coordinate on $X_0(6)$) we may write $\Psi_D=t^{\pm1}\cdot u(t)$ with $u$
 holomorphic and nonvanishing, by \eqref{eq:exact-divisor}.  A unitary
 multiplier does not affect $d\log$, and
 \[
  d\log\Psi_D=\pm\frac{dt}{t}+d\log u,
 \]
 so the residue of $\omega_D$ at $P_Q$ is $\pm1/\sqrt{-D}$, independently of
 the coordinate.  Since $D$ is odd, $\sqrt{-D}$ is a $2$-adic unit.
\end{proof}

\subsection{The partition Lambert series}

The preceding discussion describes $F_D$ geometrically.  We next compute
the same form at the cusp $\infty$.  This calculation is the promised bridge
back to the partition function.  Namely, after reduction modulo $2$, the exponents
of the Borcherds product become the partition values in \cref{thm:main}.

\begin{proposition}\label{prop:lambert}
 For square-free $D\equiv23\pmod{24}$, we have
 \begin{equation}\label{eq:lambert}
  F_D(q)\equiv
  \sum_{\substack{m\geq1\\(m,6)=1}}
  p\!\left(\frac{Dm^2+1}{24}\right)
  \sum_{\substack{n\geq1\\(n,D)=1}}q^{mn}
  \pmod2.
 \end{equation}
\end{proposition}

\begin{proof}
 Write $\chi=\leg{-D}{\cdot}$, the odd primitive quadratic character of
 conductor $D$.  Logarithmic differentiation of a single factor of
 \eqref{eq:PD} gives
 \[
  \Theta\log\bigl(1-e(-b/D)q^m\bigr)
  =-m\sum_{k\geq1}e(-bk/D)\,q^{mk},
 \]
 and summing against $\chi(b)$ produces the Gauss sums
 \[
  \sum_{b\bmod D}\chi(b)\,e(-bk/D)=\chi(-k)\tau(\chi)=-\chi(k)\sqrt{-D},
 \]
 valid for all $k$ because $\chi$ is primitive (see, for example,
 \cite[Chapters~2 and~9]{Davenport}); here
 we choose the square root $\sqrt{-D}$ so that $\tau(\chi)=\sqrt{-D}$.  Dividing by $\sqrt{-D}$,
 logarithmically differentiating \eqref{eq:PsiD} therefore gives the exact
 identity
 \begin{equation}\label{eq:exact-FD-expansion}
  F_D(q)=\sum_{m\geq1}mC(\overline m;Dm^2)
  \sum_{n\geq1}\leg{-D}{n}q^{mn}.
 \end{equation}

 We reduce \eqref{eq:exact-FD-expansion} modulo $2$ term by term.  If
 $(m,6)=1$, then $\overline m\in\{1,5,7,11\}$ and, by \eqref{eq:weil}, the
 coefficient $C(\overline m;Dm^2)$ is $\pm$ the coefficient of $q^{Dm^2}$
 in $q^{-1}f(q^{24})$, namely $\pm a_f\bigl((Dm^2+1)/24\bigr)$ where
 $f(q)=\sum a_f(n)q^n$.  Since $m$ is odd, \cref{lem:durfee} gives
 \[
  mC(\overline m;Dm^2)\equiv
  p\!\left(\frac{Dm^2+1}{24}\right)\pmod2.
 \]
 The remaining components contribute $0$ modulo $2$.  If $m$ is even, then
 $mC(\overline m;Dm^2)$ is even.  If $m$ is odd and $3\mid m$, then
 $\overline m\in\{3,9\}$ and $C(\overline m;\,\cdot\,)=0$ by
 \eqref{eq:weil}.
 Finally, we have
 \[
  \leg{-D}{n}\equiv
  \begin{cases}
   1\pmod2,&(n,D)=1,\\
   0\pmod2,&(n,D)>1.
  \end{cases}
 \]
 Substitution in \eqref{eq:exact-FD-expansion} proves
 \eqref{eq:lambert}.  
\end{proof}

\subsection{Propagation of one occurrence}

The geometric argument will produce one value of each parity.  To obtain
infinitely many, and to retain the stated bounds for the first occurrences,
we use the following arithmetic theorem
\cite[Theorem~1.2]{Ono2010}.

\begin{theorem}[Ono]\label{thm:propagation}
 Let $D>1$ be square-free with $D\equiv23\pmod{24}$.  Then the following are true.
 \begin{enumerate}
 \item If one admissible $m$ gives an even value of
 $p((Dm^2+1)/24)$, then infinitely many do.  The least such $m$ is at
 most the right-hand side of the second inequality in
 \eqref{eq:first-bounds}.
 \item If one admissible $m$ gives an odd value, then infinitely many do.
 The least such $m$ is at most $12h(-D)+2$.
 \end{enumerate}
\end{theorem}

\section{The ordinary fiber at \texorpdfstring{$2$}{2} and separation of CM points}
\label{sec:geometry}

This section supplies the arithmetic-geometric background needed to reduce
$\omega_D$ modulo $2$.  We begin with the integral model and the two
ordinary components, then explain the local orientation of a CM level
structure and the canonical-lift argument that separates its reductions.
The final subsection records the local calculation showing that a simple
pole with unit residue survives in the special fiber.  Standard references
for the integral models are Deligne--Rapoport \cite[Chapter~VI]{DR} and
Katz--Mazur \cite{KM}.

Let $\cX$ be the Deligne--Rapoport model of $X_0(6)$ over $\Z_{(2)}$.
A geometric noncuspidal point of $X_0(6)$ represents an elliptic curve $E$
together with a cyclic finite subgroup scheme $C_6\subset E$ of order $6$.
The unique primary parts of this subgroup will be denoted
\[
 C_2\subset C_6 \qquad\text{and}\qquad C_3\subset C_6.
\]
Because $3$ is invertible over $\Z_{(2)}$, the group scheme $C_3$ is
finite \'etale near the special fiber.  The geometry at $2$ is therefore
controlled by $C_2$.

Recall that an elliptic curve $\overline E$ over a field of characteristic
$2$ is \emph{ordinary} when its $2$-divisible group has a connected part of
height $1$ and an \'etale quotient of height $1$.  Over an algebraic closure
there is a connected--\'etale sequence
\begin{equation}\label{eq:connected-etale}
 0\longrightarrow \mu_{2^\infty}
 \longrightarrow \overline E[2^\infty]
 \longrightarrow \Q_2/\Z_2\longrightarrow0.
\end{equation}
In particular, $\overline E[2]$ has a unique connected subgroup of order $2$,
namely the kernel of Frobenius, and an \'etale quotient of order $2$.  A $\Gamma_0(2)$ level
structure on the ordinary locus chooses one of these two types.  This
connected--\'etale dichotomy is the source of the two components below.
For background on ordinary elliptic curves, finite flat level structures,
and the ordinary moduli problem, see
\cite[Chapters~12--13]{KM}.

Because $-D\equiv1\pmod8$, the prime $2$ splits in
$K=\Q(\sqrt{-D})$, and the CM points of discriminant $-D$ have ordinary
reduction; see, for example, \cite[Section~14]{Cox}.  The picture to keep in mind is the
following.  The special fiber of $\cX$ has two components crossing at the
supersingular points.  The
cusp $\infty$ and its $q$-expansion lie on the component where $C_2$ is
connected.  We will prove that all CM poles of $\omega_D$ land on this same
component (\cref{lem:orientation}), that no two collide there
(\cref{lem:separation}), and that every pole survives reduction
(\cref{lem:pole-survival,prop:reduced-poles}).

\subsection{The two ordinary components}

We first locate the component on which the $q$-expansion at $\infty$ lives.
The following standard description also fixes the meaning of the notation
$C_\infty$ used throughout the rest of the argument.

\begin{lemma}\label{lem:components}
 Away from the supersingular points, the normalization of the special fiber
 of $\cX$ is the disjoint union of two copies of the ordinary locus of
 $X_0(3)_{\F_2}$.  They are distinguished by whether the order-$2$ subgroup
 is connected or \'etale.  The cusp $\infty$ lies on the connected-subgroup
 component, denoted $C_\infty$.
\end{lemma}

\begin{proof}
 This is the Deligne--Rapoport description of $X_0(2N)$ at $2$, with
 $N=3$; see, for example, \cite[Chapter~VI]{DR}.  On the ordinary locus, forgetting the
 order-$2$ subgroup leaves a point of $X_0(3)$.  Conversely, the
 connected--\'etale sequence \eqref{eq:connected-etale} gives two possible geometric
 types of order $2$ subgroup: the connected kernel of Frobenius and the
 \'etale kernel of Verschiebung.  These choices define the two ordinary
 components.  After normalization, each component is therefore a
 copy of the ordinary locus of $X_0(3)_{\F_2}$.  Their images in the special
 fiber meet precisely where the elliptic curve is supersingular, which is
 also where the connected--\'etale description breaks down.

 At $\infty$, the universal generalized elliptic curve is the Tate curve
 with cyclic level subgroup $\mu_6$.  Its order $2$ part is $\mu_2$, which
 is connected in characteristic $2$.  Therefore, $\infty$ lies on the connected
 component.
\end{proof}

\subsection{The common \texorpdfstring{$2$}{2}-orientation of the two packets}

We next explain what an orientation means in this setting.  Since $2$ splits
in $K$, we write
\[
 2\cO_K=\mathfrak p_2\overline{\mathfrak p}_2,
 \qquad
 \cO_K\otimes\Z_2\simeq
 \cO_{K,\mathfrak p_2}\times
 \cO_{K,\overline{\mathfrak p}_2}
 \simeq\Z_2\times\Z_2.
\]
If $E$ has CM by $\cO_K$, this factorization decomposes its $2$-divisible
group into the summands attached to the two primes above $2$.  At a chosen
ordinary place, exactly one summand is connected and the other is \'etale.

An order-$2$ CM-stable subgroup is consequently either
$E[\mathfrak p_2]$ or $E[\overline{\mathfrak p}_2]$.  Specifying which one
occurs is the \emph{$2$-orientation} of the corresponding point of
$X_0(6)$.  This is the CM form of the connected--\'etale distinction in
\eqref{eq:connected-etale}.  See \cite[Section~14]{Cox} for the reduction of
CM elliptic curves and \cite{GrossX0N} for the oriented Heegner-point
description on $X_0(N)$.

The next lemma verifies that the two packets in the Borcherds divisor make
the same choice of prime above $2$.  It is this common orientation that puts
all their reductions on one component.

\begin{lemma}\label{lem:orientation}
 There is a fixed prime $\mathfrak p_2\mid2$ of $K$ such that the order $2$
 subgroup of every Heegner point in both packets of
 \eqref{eq:exact-divisor} is $E[\mathfrak p_2]$.
\end{lemma}

\begin{proof}
 Let $Q=[a,b,c]$ lie in $\cQ_{D,1}$ or $\cQ_{D,5}$, and write
 \[
  \tau_Q=\frac{-b+\sqrt{-D}}{2a},\qquad
  E_Q=\C/(\Z\tau_Q+\Z).
 \]
 In the standard $\Gamma_0(6)$ description, the order $2$ part of the
 cyclic level subgroup is generated by $1/2$.  Put
 \[
  \alpha_Q=a\tau_Q=\frac{-b+\sqrt{-D}}2\in\cO_K.
 \]
 The relation $a\tau_Q^2+b\tau_Q+c=0$ shows that multiplication by
 $\alpha_Q$ preserves $\Z\tau_Q+\Z$, so $\alpha_Q$ is an endomorphism of
 $E_Q$.
 Since $2\mid a$, we have
 \[
 \alpha_Q\left(\frac12\right)=\frac a2\tau_Q
 \in\Z\tau_Q+\Z.
 \]
 Put $w=(1+\sqrt{-D})/2$, so
 $\alpha_Q=w-(b+1)/2$.  In both packets $b\equiv1\pmod4$, whence as ideals of $\cO_K$,
 \[
  (2,\alpha_Q)=(2,w-1).
 \]
 The minimal polynomial of $w$ reduces to $X(X-1)$ modulo $2$, so this is a
 prime ideal of norm $2$, independent of $Q$.  We denote it by
 $\mathfrak p_2$.
 The subgroup $E_Q[\mathfrak p_2]$ has order $2$, and hence the order-$2$ level subgroup equals
 $E_Q[\mathfrak p_2]$.

The same calculation at $3$ gives
\[
 \left\langle\frac13\right\rangle
 =E_Q[(3,\alpha_Q)].
\]
Since $-D\equiv1\pmod3$, the prime $3$ splits in $K$. Here $(3,\alpha_Q)$ is one of the two primes of $\cO_K$ above $3$ and
may depend on the packet.  No common orientation at $3$ is needed.  Indeed, in the
argument below, the common orientation is required only for the order-$2$
subgroups, in order to place all the CM points on the same ordinary
component.  For the order $3$ subgroup, we use only that it is CM-stable
and that, being prime to $2$, it lifts uniquely.  Therefore, the full
$\Gamma_0(6)$ level structure is CM-stable.
\end{proof}

Choose a prime $\mathfrak{P}$ of the Hilbert class field above $\mathfrak p_2$.  At this
place, the $\mathfrak p_2$-divisible summand of $E[2^\infty]$ specializes to
the connected part of the ordinary $2$-divisible group, while the
$\overline{\mathfrak p}_2$-divisible summand specializes to the \'etale part.
This is the classical Deuring description of ordinary reduction for CM
elliptic curves; see, for example, \cite[Section~14]{Cox} and
\cite{GrossX0N}.
Consequently, $E[\mathfrak p_2]$ specializes to the connected order-$2$
subgroup, and every point of the divisor \eqref{eq:exact-divisor} lies on
$C_\infty$ after reduction.

\subsection{Canonical lifts and injectivity of specialization}

We recall the part of Serre--Tate theory used in the separation argument. 
Let $k$ be a finite field of characteristic $2$, and let
$R:=W(k)$ be its ring of Witt vectors.  Thus, $R$ is the complete
unramified discrete valuation ring of characteristic $0$ with residue field
$k$ and uniformizer $2$.  Let $\overline E/k$ be an ordinary elliptic curve.
The Serre--Tate theorem identifies
deformations of $\overline E$ over complete local $R$-algebras with residue
field $k$ with
deformations of its $2$-divisible group.  In view of
\eqref{eq:connected-etale}, such a deformation is an extension of its
\'etale part by its connected part.  More precisely, Serre--Tate theory
encodes a deformation over an algebra $A$ by a bilinear pairing between the
\'etale Tate modules of $\overline E$ and its dual, with values in
$\widehat{\mathbb G}_m(A)$.  After choosing generators, this pairing is
measured by the parameter
$q_{\mathrm{ST}}\in\widehat{\mathbb G}_m,$
called the \emph{Serre--Tate parameter}.  The identity
$q_{\mathrm{ST}}=1$ is the split extension and corresponds to the
\emph{canonical lift} of $\overline E$ (for example, see \cite[Theorems~1.2.1 and~2.1]{KatzST}).

The essential rigidity comes from endomorphisms.  If an endomorphism
$\alpha$ of $\overline E$ acts on the connected and \'etale factors of
$\overline E[2^\infty]$ by scalars
$\alpha_{\mathrm{conn}}$ and $\alpha_{\mathrm{et}}$, respectively, then
the lifting criterion in \cite[Theorem~2.1(4)]{KatzST} says that the pairing
must be equivariant for $\alpha$ and its dual $\alpha^\vee$.  For the
principally polarized CM elliptic curves below, $\alpha^\vee$ is complex
conjugation of $\alpha$, and the equivariance condition is
\begin{equation}\label{eq:ST-lifting-condition}
 q_{\mathrm{ST}}^{\alpha_{\mathrm{conn}}-
                         \alpha_{\mathrm{et}}}=1.
\end{equation}
Thus, if the difference of the two scalars is a $2$-adic unit, exponentiation
by that difference is an automorphism of $\widehat{\mathbb G}_m$ and
\eqref{eq:ST-lifting-condition} forces $q_{\mathrm{ST}}=1$.  In words, the
presence of one endomorphism whose two local actions differ by a unit forces
the ordinary deformation to be the canonical lift.

We will also use the compatibility of this lift with the level structure.
A prime-to-$2$ subgroup such as $\overline C_3$ is finite \'etale and lifts
uniquely over a complete local base.  The connected subgroup
$\overline C_2$ lifts inside the connected part of the canonical
$2$-divisible group.  It is the canonical finite flat subgroup of order $2$.
Canonical lifts with prime-to-$p$ level structures, their canonical
subgroups, and the relevant uniqueness are treated in
\cite[Sections~5.1, 5.4, and Theorem~6.1]{BorgerGurney}.

We can now prove the separation statement needed later.  The point is not
merely that the underlying CM elliptic curves have ordinary reduction: the
canonical-lift argument must also recover both parts of the $\Gamma_0(6)$
level structure.

\begin{lemma}[CM separation]\label{lem:separation}
 After choosing the place above $\mathfrak p_2$ as above, distinct points in
 the two Heegner packets of \eqref{eq:exact-divisor} have distinct geometric
 reductions on $C_\infty$.
\end{lemma}

\begin{proof}
 Individual maximal-order Heegner points with their CM-stable level
 structures are defined over the Hilbert class field $H$ of $K$; see, for
 example, \cite{GrossX0N}.  Since $H/K$ is unramified and
 $K_{\mathfrak p_2}=\Q_2$, their local fields at the chosen place are
 unramified over $\Q_2$.

 Suppose two packet points $x_1$ and $x_2$ have the same geometric reduction
 $\overline x=(\overline E,\overline C_2,\overline C_3)$.  After a finite
 unramified extension of the common $2$-adic field, we may identify both
 special fibers with $\overline x$; both lifts and all their CM endomorphisms
 are then defined.  Specialization of endomorphisms gives an injection $\cO_K \hookrightarrow \End(\overline E)$.  Since $\overline E$ is ordinary,
 $\End(\overline E)\otimes\Q$ is an imaginary quadratic field.  It must be
 $K$, and maximality of $\cO_K$ gives $\End(\overline E)=\cO_K$.  Therefore, the two
 specialized CM actions agree up to complex conjugation.

 Every endomorphism in $\cO_K$ lifts on each $E_i$.  Since $2$ splits,
 \[
  \cO_K\otimes\Z_2\simeq\Z_2\times\Z_2,
 \]
 the two factors corresponding to the actions of $\cO_K$ on the connected
 and \'etale parts of $\overline E[2^\infty]$.  In particular, complex conjugation on
 $\cO_K$ interchanges them.  For $\alpha\in\cO_K$, write
 $\alpha_{\mathrm{conn}}$ and $\alpha_{\mathrm{et}}$ for its two components.
 Under the principal polarization, the dual action of $\alpha$ is the action
 of $\overline\alpha$.  Hence the lifting condition is precisely
 \eqref{eq:ST-lifting-condition}.
 Take $\alpha=w=(1+\sqrt{-D})/2$.  Then
 $\alpha_{\mathrm{conn}}-\alpha_{\mathrm{et}}=\pm\sqrt{-D}$ is a
 $2$-adic unit, because $D$ is odd. The sign accounts for the two possible
 identifications of the CM action with $\End(\overline E)$.
 Exponentiation by a $2$-adic unit is an automorphism of
 $\widehat{\mathbb G}_m$, so $q_{\mathrm{ST}}=1$.  Therefore, both $E_1$ and
 $E_2$ are the canonical lift of $\overline E$.

 In particular, there is an isomorphism $E_1\simeq E_2$ reducing to the
 identity of $\overline E$.  The subgroup $\overline C_3$ lifts uniquely to
 either side because $3$ is invertible and
 $E[3]$ is finite \'etale over the deformation ring, so the isomorphism
 carries $C_3^{(1)}$ to $C_3^{(2)}$.  The subgroup
 $\overline C_2$ lifts as the unique finite flat order $2$ subgroup scheme with
 connected special fiber, namely the order $2$ subgroup of the connected
 part of the canonical $2$-divisible group. Then by \cref{lem:orientation}, this
 is $E[\mathfrak p_2]$.  Therefore, the isomorphism also carries
 $C_2^{(1)}$ to $C_2^{(2)}$.  So it respects both parts of the
 $\Gamma_0(6)$ structure, and $x_1=x_2$, a contradiction.
\end{proof}

\subsection{Survival of poles with unit residues}

We finally explain what it means to reduce a meromorphic differential.  Let
$R$ be an unramified discrete valuation ring over $\Z_2$ and let $\mathcal Y$
be a regular relative curve over $R$.  A point of the generic fiber may
extend to a \emph{horizontal section} of $\mathcal Y$, while an irreducible
component of the special fiber is a \emph{vertical divisor}.  Before a
rational differential on the generic fiber can be reduced on such a
component, one must know its order along that vertical divisor.

At the generic point of a smooth component $C$, the local ring of
$\mathcal Y$ is a discrete valuation ring with uniformizer $2$.  The
\emph{vertical order} of a differential $\omega$ along $C$ is the valuation
of $\omega/dt$ there, where $t$ reduces to a separating local parameter on
$C$.  Equivalently, for an expansion
\[
 \omega=\left(\sum_n a_nt^n\right)dt
\]
at a smooth point of $C$, it is the $2$-adic Gauss valuation
$\min_n v_2(a_n)$.  Nonnegative vertical order permits reduction.  In particular, order zero
gives a nonzero rational differential on $C$.  The next lemma combines two
simple but crucial observations: an integral expansion gives nonnegative
vertical order, while a unit residue at a separated horizontal pole prevents
the order from being positive.  This allows us to reduce $\omega_D$ without
trying to reduce the Borcherds product $\Psi_D$ itself.

\begin{lemma}[Pole survival]\label{lem:pole-survival}
 Let $R$ be a complete discrete valuation ring, unramified over $\Z_2$ (so
 $2$ is a uniformizer), and
 let $C$ be a smooth irreducible component of the special fiber of a regular
 relative curve over $R$.  Let $\omega$ be a rational differential on the
 generic fiber satisfying the following conditions:
 \begin{enumerate}
 \item Its expansion at an $R$-point of $C$ is integral, where the reduction
 of this $R$-point is a smooth point of the special fiber lying on $C$ alone
 and meeting no polar section of $\omega$.
 \item It has a simple pole along a horizontal section, defined over $R$,
 meeting $C$ at a smooth point of the special fiber, with unit residue.
 \item No other pole section meets the same point of $C$.
 \end{enumerate}
 Then $\omega$ has vertical order zero along $C$, and its reduction on $C$
 has a simple pole with nonzero residue.
\end{lemma}

\begin{proof}
 At the reduction of the $R$-point in condition (1), the completed local
 ring of the total space is $R[[t]]$, with $C$ cut out by $2$, because the
 point is a smooth point of the reduced special fiber lying on the single
 component $C$.  Since no polar section of $\omega$ passes through this
 point, the only prime divisor of $\operatorname{Spec}R[[t]]$ along which
 $\omega/dt$ can have a pole is $C=(2)$; because $R[[t]]$ is a regular
 local ring, and hence a unique factorization domain, it follows that
 $\omega/dt$ lies in $R[[t]][1/2]$, and its expansion there computes its
 Gauss valuation.
 Integrality of the expansion therefore implies that the vertical order of
 $\omega$ along $C$ is nonnegative.

 At the intersection with the horizontal pole
 section, the completed local ring is again $R[[t]]$, with the section defined
 by $t=0$.  Since no other horizontal pole passes through the point, the pole
 along the section is simple.  Because the vertical order is nonnegative,
 $t\omega/dt$ belongs to $R[[t]]$.  Its value
 at $t=0$ is the residue $u\in R^\times$. Hence $t\omega/dt=u+t g(t)$ for some $g(t)\in R[[t]]$, and
 \[
  \omega=\left(\frac{u}{t}+g(t)\right)dt,
  \qquad u\in R^\times,\quad g(t)\in R[[t]].
 \]
 If the vertical order were positive, every coefficient, including the
 residue $u$, would be divisible by the uniformizer.  Hence the vertical
 order is zero.  Reducing the displayed expansion gives
 $\overline u\,dt/t$ plus a regular differential, and
 $\overline u\neq0$.
\end{proof}

Combining the global integrality at the cusp with the local information at
the separated CM sections gives the form of pole survival that will be used
in the main proof.

\begin{proposition}\label{prop:reduced-poles}
 The differential $\omega_D$ has vertical order zero along $C_\infty$.
 Its reduction $\overline\omega_D$ is a nonzero rational differential on
 $C_\infty$ with a simple pole at the reduction of every point in
 \eqref{eq:exact-divisor}.
\end{proposition}

\begin{proof}
 Work over $R=\cO_{H_\mathfrak{P}}$, the ring of integers of the completion of the
 Hilbert class field at the chosen prime $\mathfrak{P}$ above $\mathfrak p_2$.  It is
 unramified over $\Z_2$, as noted in the proof of \cref{lem:separation}.
 The coarse model $\cX_R$ can fail to be regular at the supersingular
 points of its special fiber, where the reduced curve has extra
 automorphisms, so we apply \cref{lem:pole-survival} on the open complement
 $\cX_R^{\mathrm{ord}}$ of these points.  Away from the supersingular
 points the special fiber is smooth, so $\cX_R^{\mathrm{ord}}$ is smooth
 over $R$; in particular, it is a regular relative curve, and it contains
 the cusp $\infty$ and, by ordinary reduction, every CM section considered
 below.
 The cusp $\infty$ is an $R$-point of $\cX_R$ reducing to a smooth point of
 the special fiber on $C_\infty$ alone, and it meets no polar section of
 $\omega_D$. Singular moduli are integral, and the splitting of $2$ in $K$
 gives ordinary reduction, so the CM sections pass through noncuspidal smooth
 points of $C_\infty$.
 Since $F_D=\sum_{n\geq1}a_nq^n$ has integral coefficients and no constant
 term,
 \[
  \omega_D=\biggl(\sum_{n\geq1}a_nq^{n-1}\biggr)dq
 \]
 has integral
 expansion at $\infty$, supplying condition (1) of
 \cref{lem:pole-survival}.  Each CM point of \eqref{eq:exact-divisor}
 extends, by properness, to an $R$-point of $\cX_R$. These are the
 horizontal pole sections, with unit residues by
 \cref{cor:unit-residues}.  They are pairwise
 disjoint on $C_\infty$ by \cref{lem:separation}, giving conditions (2)
 and (3).  Apply \cref{lem:pole-survival} on $\cX_R^{\mathrm{ord}}$ at
 every CM section.
\end{proof}

\section{Proof of Theorem~\ref{thm:main}}\label{sec:main-proof}

We now bring together the two descriptions of $F_D$.  The CM divisor will
rule out an identically zero reduction, while a holomorphic Eisenstein
series will rule out the possibility that every relevant coefficient is
one.  The propagation theorem then converts these two initial occurrences
into the full conclusion of \cref{thm:main}.
For brevity, write
\[
 A_D(m):=p\!\left(\frac{Dm^2+1}{24}\right)
 \qquad ((m,6)=1).
\]

\subsection{An odd value}

The first half of the argument uses only the Lambert series and the
surviving CM poles.
Suppose first that every $A_D(m)$ is even.  By \cref{prop:lambert}, the
formal expansion of $\overline\omega_D$ at $\infty$ is identically zero.
A rational differential on the irreducible curve $C_\infty$ is determined by
its expansion on any nonempty formal neighborhood.  Therefore,
$\overline\omega_D$ would vanish identically.  This contradicts
\cref{prop:reduced-poles}.  Therefore at least one $A_D(m)$ is odd.

\subsection{An even value}

To rule out constant odd parity, we construct a holomorphic modular form
whose reduction has exactly the Lambert series that would arise if every
$A_D(m)$ were odd.
Let $E_2(z)=1-24\sum_{n\geq1}\sigma_1(n)q^n$ and define
\begin{equation}\label{eq:E}
 E(z):=\frac{(E_2(z)-3E_2(3z))-2(E_2(z)-2E_2(2z))}{24}.
\end{equation}
Each difference $E_2(z)-rE_2(rz)$ is a holomorphic weight $2$ modular form
on $\Gamma_0(r)$, so $E\in M_2(\Gamma_0(6))$.  It has zero constant term,
and, with the convention that $\sigma_1(x)=0$ for nonintegral $x$, its
coefficient of $q^n$ is the integer
\[
 \sigma_1(n)-4\sigma_1(n/2)+3\sigma_1(n/3).
\]
Write $n=2^a3^bu$ with $(u,6)=1$.  Modulo $2$, this coefficient is
$\sigma_1(n)+\sigma_1(n/3)$.  Now $\sigma_1(2^a)$ is odd.  If $b\geq1$, then
\[
 \sigma_1(3^b)+\sigma_1(3^{b-1})\equiv3^b\equiv1\pmod2,
\]
while for $b=0$ the second term is absent.  Finally,
$\sigma_1(u)\equiv d(u)\pmod2$, since every divisor of $u$ is odd.  Therefore, the
coefficient is congruent to the number $d(u)$ of divisors of $u$.
Consequently, we have
\begin{equation}\label{eq:E-mod2}
 E(z)\equiv
 \sum_{\substack{m\geq1\\(m,6)=1}}
 \sum_{n\geq1}q^{mn}\pmod2.
\end{equation}
We define
\begin{equation}\label{eq:ED}
 E_D(z):=\sum_{\delta\mid D}E(\delta z)\in M_2(\Gamma_0(6D)).
\end{equation}
In the term coming from $E(\delta z)$, the second index in
\eqref{eq:E-mod2} is multiplied by $\delta$.  Thus, for fixed $m$ and fixed
second index $n$, the coefficient is multiplied by
$\sum_{\delta\mid(D,n)}1$, the number of divisors of $(D,n)$.  Since $D$ is
square-free, this count is a power of $2$, and so it is odd exactly when
$(D,n)=1$.  Therefore, we have
\begin{equation}\label{eq:ED-mod2}
 E_D(z)\equiv
 \sum_{\substack{m\geq1\\(m,6)=1}}
 \sum_{\substack{n\geq1\\(n,D)=1}}q^{mn}\pmod2.
\end{equation}
This is \cite[equation~(3.8)]{Ono2010}.

Suppose that every $A_D(m)$ is odd.  Equations
\eqref{eq:lambert} and \eqref{eq:ED-mod2} then imply
\begin{equation}\label{eq:F-equals-E}
 F_D(q)\equiv E_D(q)\pmod2.
\end{equation}

Let $\pi:X_0(6D)\to X_0(6)$ forget the cyclic subgroup of order $D$.
On the ordinary noncuspidal locus at $2$, the forgetful map associated
with the odd-order $D$-level is finite \'etale, and it does not change the connected--\'etale
dichotomy for the order $2$ subgroup.  The component $C'_\infty$ of the
normalized special fiber of $X_0(6D)$ containing $\infty$ maps surjectively to
$C_\infty$.  Equivalently, its normalization is the component corresponding to the
connected choice of the order-$2$ subgroup, and is identified on the
ordinary locus with \(X_0(3D)_{\F_2}\). The cusp $\infty$ has width $1$ on both curves and
$\pi$ is formally the identity in the parameter $q$ there, so the expansion
of $\pi^*\omega_D$ at $\infty$ is again $F_D(q)\,dq/q$.  Since $\pi$ is
generically \'etale on $C'_\infty$, the vertical order of $\pi^*\omega_D$
along $C'_\infty$ equals the vertical order of $\omega_D$ along $C_\infty$,
which is zero by \cref{prop:reduced-poles}, and the reduction of
$\pi^*\omega_D$ on $C'_\infty$ is the pullback of $\overline\omega_D$.
Choose a point $y\in C'_\infty$ above one of the CM poles of
$\overline\omega_D$.  Since \(D\) is odd, the additional \(D\)-level is prime to the residue
characteristic, and the underlying CM point is nonelliptic because \(-D<-4\).  Hence the forgetful map is \'etale at \(y\).
Therefore, $\pi^*\overline\omega_D$ has a simple pole at $y$.

On the other hand, the reduction on $C'_\infty$ of the differential
$E_D(q)\cdot \frac{dq}{q}$
is regular at $y$.  Indeed, the integral expansion of $E_D$ at $\infty$
gives nonnegative vertical order along $C'_\infty$, while the differential associated with the holomorphic weight-$2$
modular form $E_D$ has no horizontal poles at noncuspidal points. Indeed, $X_0(6D)$ has no elliptic points: the standard criterion gives no elliptic points of order $2$ because $3\mid6D$ and $3\equiv3\pmod4$, and none of order $3$ because $2\mid6D$ and $2\equiv2\pmod3$. Thus, the reduction is either
zero or a rational differential regular at $y$, and in either case its formal
expansion at $\infty$ is the reduction of $E_D(q)\,dq/q$ modulo $2$.
Congruence
\eqref{eq:F-equals-E} says that the two reduced differentials have the same
formal expansion at $\infty$.  Hence they agree as rational differentials on
the irreducible curve $C'_\infty$.  This is impossible because one has a pole
at $y$ and the other is regular there.  Therefore, at least one $A_D(m)$ is even.

We have proved one occurrence of each parity.  Applying
\cref{thm:propagation} gives 
 infinitely many occurrences of both, as well as the claimed bounds for
$m_{\mathrm{odd}}$ and $m_{\mathrm{even}}$.

\section{Sketch of the proof of Theorem~\ref{thm:general}}\label{sec:general-proof}

Let $L/\Q_2$ be the common unramified field in the hypotheses and let $R$ be
its ring of integers.  Choose the CM place above the fixed prime
$\mathfrak p_2$ in \emph{(G2)}.  At this place $E[\mathfrak p_2]$ is
connected, so every CM pole lies on the ordinary component $C_\infty$
containing the Tate cusp.  The prime-to-$2$ level structure does not affect
this choice of component.

Distinct poles remain distinct after reduction.  Indeed, suppose that two
of them reduce to the same ordinary enhanced elliptic curve $\overline E$.
As in the proof of \cref{lem:separation}, specialization of the
maximal-order CM action and ordinarity imply \(\End(\overline E)=\cO_K\), and the CM
actions lift on both deformations.  Since
$\cO_K\otimes\Z_2\simeq\Z_2\times\Z_2$, one may choose an endomorphism whose
actions on the connected and \'etale parts differ by a unit.  The
Serre--Tate lifting criterion then forces both deformations to be the
canonical lift of $\overline E$.  Their prime-to-$2$ subgroup schemes lift
uniquely, while their order-$2$ subgroups are the canonical connected lifts prescribed by the common \(2\)-orientation in \emph{(G2)}.  The two
enhanced CM points would therefore coincide already in
characteristic zero, a contradiction.  This is precisely where the
maximal-order hypothesis and the common $2$-orientation are used.

By \emph{(G1)}, the integral expansion of $\omega_D$ at $\infty$ gives
nonnegative vertical order along $C_\infty$.  Each separated CM section has,
by \emph{(G2)}, a simple pole with unit residue.  The local calculation in
\cref{lem:pole-survival}, applied on the ordinary locus as in the proof of
\cref{prop:reduced-poles}, therefore forces the vertical order to be zero and
shows that the reduction of $\omega_D$ has a genuine CM pole.  If every
$B_D(m)$ were even, however, \emph{(G1)} would make its formal expansion at
$\infty$ vanish.  The reduced differential would then vanish identically on
$C_\infty$, contradicting the surviving pole.  Hence, at least one $B_D(m)$
is odd.

It remains to rule out constant odd parity.  Suppose that every $B_D(m)$ is
odd, and pull back to $X_0(ND)$.  By \emph{(G1)} and \emph{(G3)}, the
differentials $\omega_D$ and $\mathcal E_D(q)dq/q$ have the same reduced
formal expansion at $\infty$.  They must therefore agree as rational
differentials on the ordinary component above $C_\infty$.  The pullback of
$\omega_D$ retains a simple CM pole: since $2$ splits in
$K=\Q(\sqrt{-D})$, we have $D\equiv7\pmod8$, so $D>4$. The CM point is
nonelliptic, and the forgetful map is \'etale there because $D$ is odd and
coprime to $N$.  On the other hand, the integral expansion of
$\mathcal E_D$ gives nonnegative vertical order, and its holomorphy excludes
a horizontal pole at that noncuspidal point.  Its reduction is therefore
regular there, a contradiction.  Therefore, at least one $B_D(m)$ is even, which
proves nonconstancy.

\appendix

\section{Proof of \texorpdfstring{\Cref{cor:sqrt}}{the Corollary}}
\label{sec:appendix}

\begin{center}
\textsc{by Qi-Yang Zheng}\\
\medskip
\texttt{zhengqy29@mail2.sysu.edu.cn}
\end{center}

\medskip

We first recall the statement to be proved.

\medskip
\noindent\textbf{\Cref{cor:sqrt}.}
\emph{If we let $N_{\mathrm{odd}}(X):=\#\{0\leq n\leq X: p(n)\text{ is odd}\}$,
then we have}
\[
 \liminf_{X\to\infty}\frac{N_{\mathrm{odd}}(X)}{\sqrt X}
 \ \geq\ \frac{243}{64\sqrt6\,\pi^5}.
\]
\emph{In particular, we have $N_{\mathrm{odd}}(X)\gg\sqrt X$.}

\medskip

\begin{proof}
Throughout, $Y$ is a large real parameter, and we set
\[
 \cD(Y):=\{D\in\Z_{>0}: D\leq Y,\ D\equiv23\!\!\pmod{24},\ D\text{ square-free}\},
 \qquad
 \delta:=\frac{3}{8\pi^2},\qquad \beta:=\frac{\pi}{108}.
\]
The idea is straightforward.  For each $D\in\cD(Y)$, \cref{thm:main} produces
an odd value of the partition function at the integer $n_D=(Dm_D^2+1)/24$ for
some admissible $m_D\leq12h(-D)+2$, and distinct $D$ turn out to give distinct
$n_D$.  Class numbers are small on average, so for most $D\in\cD(Y)$ the
integer $n_D$ is at most a fixed multiple of $Y^2$.  Comparing the number of
such $D$ with the length of the interval containing the $n_D$ yields the
$\sqrt X$ bound, and optimizing the class number threshold gives the constant.

We begin with the count $\#\cD(Y)=\delta Y+O(\sqrt Y)$.  Inserting
$\mu^2(n)=\sum_{d^2\mid n}\mu(d)$ and writing $n=d^2r$, we note that every
$n\equiv23\pmod{24}$ is coprime to $6$, hence so is $d$. Then we have
$d^2\equiv1\pmod{24}$, so the congruence on $n$ is equivalent to
$r\equiv23\pmod{24}$.  Therefore, we have
\[
 \#\cD(Y)
 =\sum_{\substack{d\leq\sqrt Y\\ (d,6)=1}}\mu(d)
   \left(\frac{Y}{24d^2}+O(1)\right)
 =\frac{Y}{24}\prod_{p\geq5}\left(1-\frac1{p^2}\right)+O(\sqrt Y)
 =\delta Y+O(\sqrt Y),
\]
since $\prod_{p\geq5}(1-p^{-2})
=\zeta(2)^{-1}(1-2^{-2})^{-1}(1-3^{-2})^{-1}=9/\pi^2$.

Next we show that
\begin{equation}\label{eq:class-average}
 \sum_{D\in\cD(Y)}h(-D)\ \leq\ (\beta+o(1))\,Y^{3/2}.
\end{equation}
For $D\in\cD(Y)$ the discriminant $-D$ is fundamental, so $h(-D)$ counts
classes of primitive positive definite binary quadratic forms of discriminant
$-D$.  Every such class contains a representative $[a,b,c]$ with
$|b|\leq a\leq c$, and exactly one such representative once one imposes in
addition that $b\geq0$ whenever $|b|=a$ or $a=c$ (for example, see
\cite[Theorem~2.8]{Cox}).  We need only an upper bound, so it is enough that
every class be counted at least once: dropping the boundary conditions, as
well as primitivity and the square-freeness of $D$, can only increase the
count.  It therefore suffices to bound the number of integer triples
$(a,b,c)$ with $|b|\leq a\leq c$, $0<4ac-b^2\leq Y$, and
$4ac-b^2\equiv23\pmod{24}$.  The congruence forces $b$ to be odd; and since
$b^2\equiv1\pmod{24}$ when $(b,6)=1$ while $b^2\equiv9\pmod{24}$ when
$3\mid b$, it amounts to $ac\equiv0\pmod 6$ in the first case and
$ac\equiv2\pmod 6$ in the second.  Thus, for fixed $a$, each odd residue class
of $b$ modulo $6$ confines $c$ to a set of residue classes modulo $6$, and a
short computation shows that the total weight $R(a)$ of these classes,
counting the class-count for $ac\equiv0\pmod6$ twice (once for each of
$b\equiv1,5\pmod6$), is periodic in $a$ modulo $6$ with the values
$12,3,6,6,6,3$ for $a\equiv0,1,\ldots,5\pmod6$. In particular, $R$ has mean
$6$, so the congruence conditions hold for one sixth of the pairs $(b,c)$ on
average.

Reducedness gives $Y\geq 4ac-b^2\geq4a^2-a^2=3a^2$, so
$a\leq\sqrt{Y/3}$. For fixed $a$ and odd $b$ with $|b|\leq a$, the
conditions $c\geq a$ and $4ac-b^2\leq Y$ confine $c$ to an interval of length
$(Y+b^2-4a^2)_+/4a$, where $x_+:=\max(x,0)$.  Counting the admissible residue classes of $c$ in this
interval and summing over $b$, the number of triples with first coefficient
$a$ is at most $\tfrac{R(a)}{6}F(a)+O(Y/a+a)$, where
\[
 F(a):=\frac{1}{24a}\int_{-a}^{a}\bigl(Y+t^2-4a^2\bigr)_+\,dt .
\]
Summing over $a\leq\sqrt{Y/3}$, partial summation replaces the periodic factor
$R(a)/6$ by its mean value $1$ at the cost of $O(Y)$, while the error terms
contribute $O(Y\log Y)$.  Substituting $a=\sqrt Y\,\alpha$ and
$t=\sqrt Y\,\alpha\sigma$ with $\sigma\in[-1,1]$, the inner $\alpha$-integral
evaluates to $\tfrac23(4-\sigma^2)^{-1/2}$, and
$\int_{-1}^{1}(4-\sigma^2)^{-1/2}\,d\sigma
=[\arcsin(\sigma/2)]_{-1}^{1}=\pi/3$, whence
\[
 \int_0^{\infty}F(a)\,da
 =\frac{Y^{3/2}}{24}\int_{-1}^{1}\int_0^{(4-\sigma^2)^{-1/2}}
   \bigl(1-(4-\sigma^2)\alpha^2\bigr)\,d\alpha\,d\sigma
 =\frac{Y^{3/2}}{24}\cdot\frac23\cdot\frac{\pi}{3}
 =\frac{\pi}{108}\,Y^{3/2}.
\]
This proves \eqref{eq:class-average}.

Now fix $H>\beta/\delta$.  By \eqref{eq:class-average}, at most
$(\beta/H+o(1))Y$ of the $D\in\cD(Y)$ satisfy $h(-D)>H\sqrt Y$, and so at
least $(\delta-\beta/H+o(1))Y$ of them satisfy $h(-D)\leq H\sqrt Y$. We call
these $D$ \emph{good}.  For each good $D$, \cref{thm:main} provides an
admissible $m_D\leq12h(-D)+2\leq12H\sqrt Y+2$ for which $p(n_D)$ is odd, where
\[
 n_D:=\frac{Dm_D^2+1}{24}\ \leq\ \frac{Y\bigl(12H\sqrt Y+2\bigr)^2+1}{24}
 =6H^2Y^2+2HY^{3/2}+\frac Y6+\frac1{24}=:Z_H(Y).
\]
The map $D\mapsto n_D$ is injective: from $n_D$ we recover
$24n_D-1=Dm_D^2$, and since $D$ is square-free, $D$ is exactly the product of
the primes dividing $Dm_D^2$ to an odd power, hence is determined by $n_D$.
The good $D$ therefore produce at least $(\delta-\beta/H+o(1))Y$ distinct
integers $n\leq Z_H(Y)$ at which $p(n)$ is odd.  That is, we have
$N_{\mathrm{odd}}(Z_H(Y))\geq(\delta-\beta/H+o(1))Y$.

Finally, given $X$, let $Y_X$ be the unique real solution of $Z_H(Y_X)=X$,
which exists because $Z_H$ is continuous and strictly increasing on
$(0,\infty)$.  The leading term of $Z_H$ gives $Y_X/\sqrt X\to1/(\sqrt6\,H)$,
whence
\[
 \liminf_{X\to\infty}\frac{N_{\mathrm{odd}}(X)}{\sqrt X}
 \ \geq\ \frac{\delta-\beta/H}{\sqrt6\,H}.
\]
The right-hand side is maximized at $H=2\beta/\delta=4\pi^3/81$, which is
admissible since $2\beta/\delta>\beta/\delta$, and the maximum value is
\[
 \frac{\delta^2}{4\beta\sqrt6}
 =\frac{\bigl(3/(8\pi^2)\bigr)^2}{4\bigl(\pi/108\bigr)\sqrt6}
 =\frac{243}{64\sqrt6\,\pi^5}.
\]
This proves the corollary.
\end{proof}

\section{Lean formalization}
\label{sec:lean}
The algebraic identities required to prove \cref{thm:main} have been
formalized and verified in Lean by AxiomProver, and the details can be found at
\begin{center}
\url{https://github.com/AxiomMath/PartitionParity}.
\end{center}
The formalization covers seven identities: the congruence
$P(q)\equiv f(q)\pmod4$ of \cref{lem:durfee}, the exact expansion
\eqref{eq:exact-FD-expansion} of $F_D$ together with its mod $2$ reduction
\eqref{eq:lambert}, the residue computation of \cref{cor:unit-residues}, the
Eisenstein congruences \eqref{eq:E-mod2} and \eqref{eq:ED-mod2}, and the
congruence \eqref{eq:F-equals-E} under the hypothesis that all of the relevant
partition values are odd.  Each identity is verified as an implication from the
paper's classical inputs. Two classical facts enter as designated hypotheses
(the Durfee-square identity of \cref{lem:durfee} and the Gauss-sum evaluation
$\tau(\chi)=\sqrt{-D}$, encoded algebraically by the relation
$\tau(\chi)^2=-D$), while the modular and geometric objects enter as opaque
data subject only to the finitely many algebraic relations actually used, such
as the coefficient table \eqref{eq:weil} and the local factorization
$\Psi_D=t^{\pm1}u(t)$ at a point of the divisor.  It does not include the
analytic and geometric inputs: the construction of the Borcherds lift and its
divisor (\cref{thm:borcherds}), the characteristic $2$ geometry of
Section~\ref{sec:geometry} (canonical lifts, separation of CM points, and
survival of poles), the parity-elimination arguments of
Section~\ref{sec:main-proof} that consume these identities, and the
propagation theorem (\cref{thm:propagation}).  

Furthermore, we note that
\cref{cor:sqrt} has also been formalized in Lean by AxiomProver, and the
results can be found at
\begin{center}
\url{https://github.com/AxiomMath/PartitionZheng}.
\end{center}

\end{document}